# High-Precision, Fair University Course Scheduling During a Pandemic

**Matthew E. H. Petering[a,1]**

**Mohammad Khamechian[a]**

[a] *Department of Industrial & Manufacturing Engineering, University of Wisconsin-Milwaukee, P.O. Box 784, Milwaukee, WI 53201, USA*

E-mail addresses:

mattpete@uwm.edu (M.E.H. Petering), mr.khamechian@gmail.com (M. Khamechian).

[1] Corresponding author. Tel.: +1 414 229 3448; fax: +1 414 229 6958.

High-Precision, Fair University Course Scheduling During a Pandemic

Abstract

***Problem definition*:** Scheduling university courses is extra challenging when classroom capacities are reduced because of social distancing requirements that are implemented in response to a pandemic such as COVID-19. ***Methodology:*** In this work, we propose an expanded taxonomy of course delivery modes, present an integer program, and develop a course scheduling algorithm to enable all course sections—even the largest—to have a significant classroom learning component during a pandemic. Our approach is fair by ensuring that a certain fraction of the instruction in every course section occurs in the classroom. Unlike previous studies, we do not allow rotating attendance and instead require simultaneous attendance in which all students in a section meet in 1-5 rooms at the same time but less often than in a normal semester. These mass meetings, which create opportunities for in-person midterm exams and group activities, are scheduled at high precision across all days of the semester rather than a single, repeating week. A fast heuristic algorithm makes the schedule in an hour. ***Results*:** We consider the 1834 in-person course sections, 172 classrooms, and 96 days in the fall 2022 semester at [*UniversityXYZ*]. If average classroom capacity is reduced by 75% due to a pandemic, our approach still allows at least 25% of the instruction in every section, and more than 49% of all instruction across the entire campus, to be in the classroom. Our method also produces excellent results for regular classroom assignment. ***Managerial implications*:** An algorithm based on the principles of fairness and simultaneous attendance can significantly improve university course schedules during a pandemic and in normal times. High-precision schedules that prepare a campus for various pandemic possibilities can be created with minimal administrative effort and activated at a moment's notice before or during a semester if an outbreak occurs.

1. Introduction

In the early 2020s, learning at the world's 25,000 colleges and universities was in jeopardy as SARS-CoV-2 (i.e., COVID-19) spread across the globe. As of May 10, 2023, COVID-19 had infected more than 765 million people and killed 6.9 million (WHO, 2023). Although many experts believe the worst of the pandemic is over, the disruption it has caused to businesses, governments, healthcare systems, and the cognitive development of young people under the age of 25 will likely be felt for decades (Schady et al., 2023). Moreover, there is still the possibility of a future pandemic.

For universities, a pandemic such as COVID-19 creates a dilemma of education versus health which requires them to rethink how their courses are delivered. If education is the priority, courses should be delivered in person (i.e., taught in a classroom) to improve learning outcomes (Freeman et al., 2014). On the other hand, if health is the only consideration, courses should be delivered online (i.e., remotely) to minimize the risk of infection.

One option for addressing this dilemma is to socially distance, i.e., to require students to remain at least six feet apart when sitting in classrooms to comply with recommendations made by the U.S. Centers for Disease Control and Prevention. Although social distancing is an obvious solution, it creates a major problem: classroom capacities are reduced by about 70% (Reeves, 2020). This in turn creates a new task for the university to allocate scarce classroom space to hundreds of courses whose total demand for space exceeds the supply.

Since the onset of COVID-19, many ideas have been proposed for managing and scheduling university instruction if social distancing is required in the classroom. The most common idea is to deliver courses in hybrid mode with rotating attendance. In such a course, students are divided into cohorts with the number of students in each cohort no greater than the pandemic capacity (i.e., socially distanced capacity) of the classroom. Each time the course meets, one cohort attends in the classroom and the other cohorts attend online, and the in-person attendance rotates among the cohorts. Another idea is to change the rooms assigned to courses so more courses can be taught in person. Here, the general idea is to move large courses online so large classrooms can be used to teach medium-sized courses in person. This in turn frees up medium-sized classrooms for teaching small courses in person.

However, the above strategies have some disadvantages. A hybrid delivery mode with rotating attendance does not provide opportunities for in-person midterm exams in which all students are simultaneously in a classroom. There is also a nonuniform student experience because of mixed attendance at each lecture; some students receive the lecture in a classroom while others receive it online. The delivery of mixed-attendance lectures may also create more work and stress for instructors as they simultaneously attempt to connect with students in two different environments. Furthermore, moving large courses online for the benefit of smaller courses could be perceived as unfair. Why should the courses with the highest enrollments—the most successful courses by some metrics—be sacrificed for the sake of smaller courses? In this study, we propose an alternate framework for delivering and scheduling university courses which addresses these issues. Our framework enables all course sections—even the largest—to have a significant classroom learning component during a pandemic.

*1.1. The course scheduling process*

Before discussing our approach in detail, we first place it in the context of the overall course scheduling process. The tasks involved in scheduling a semester of courses at a major U.S. university are shown in the left column of Figure 1. First, faculty members identify the subjects they will teach. This decision does not change much from one semester to another. Next, about seven months before the start of a semester, each department decides what courses it will offer that semester. This is called term planning. Soon thereafter, each department decides the preferred meeting times of its courses. These times are chosen to maximize student and instructor satisfaction, and they are typically selected manually. About 5-6 months before the semester begins, the university registrar uses a computer program to assign a room to as many courses as possible, assuming each course meets during its preferred meeting time. Typically, some courses with preferred meeting times during peak hours (e.g., Mon-Thurs from 9 AM to 3 PM) are unable to be assigned a room, so the registrar contacts some departments and asks them to change the meeting times of those courses. Hence, there is a dotted arrow from task 4 to 3 in the figure. After a few iterations of tasks 3-4, the registrar is able to feasibly assign a room to every course. Later, students enroll in courses in waves based on seniority.

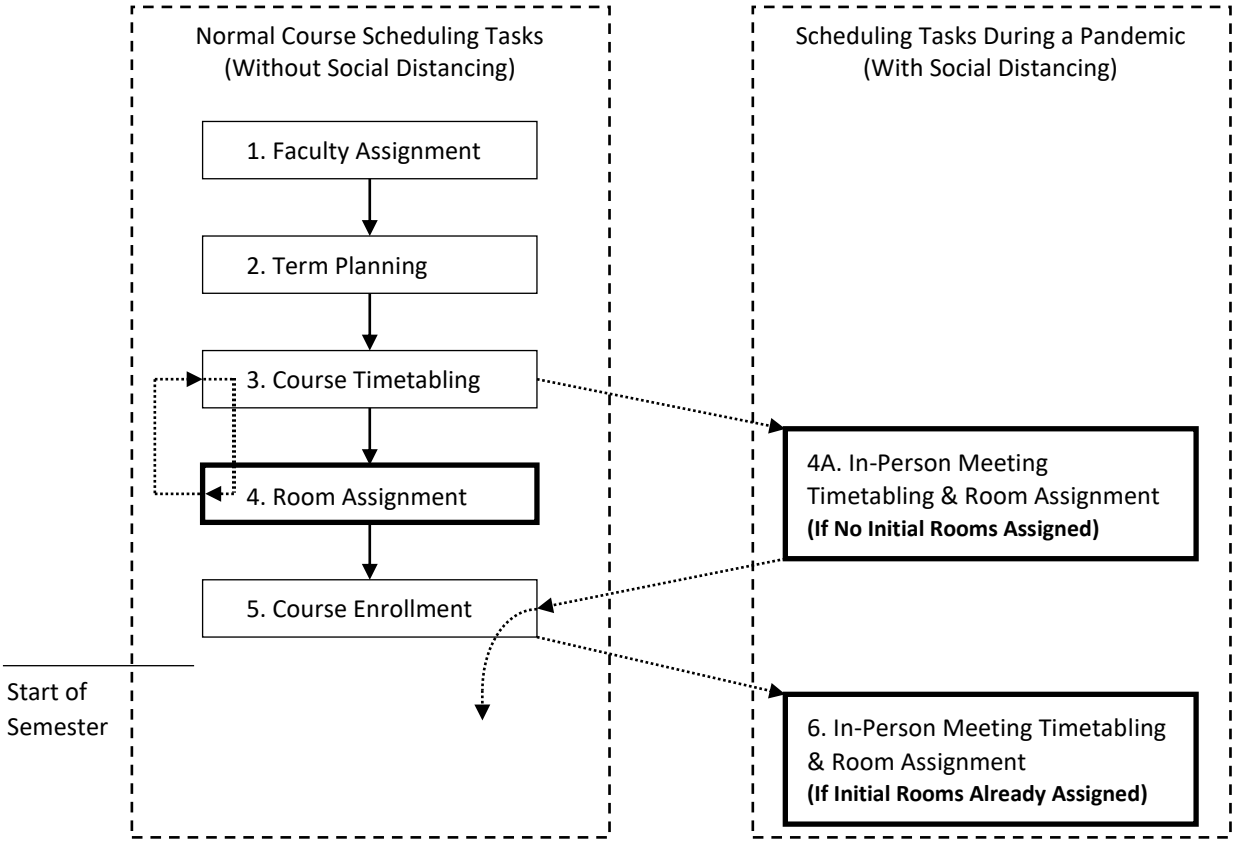

**Figure 1.** The university course scheduling process from Barnhart et al. (2022) with additions. Normal course scheduling tasks are in the left column. Tasks involved in preparing for, or responding to, a pandemic are in the right column. Dotted arrows show possible decision pathways. The tasks addressed in this study are in bold.

In this work, we propose a course scheduling approach that handles the tasks in bold in Figure 1. First, if a university wishes to prepare for, or respond to, a pandemic after tasks 1-5 are completed, our approach can decide when and where each course meets in a socially distanced manner from that point onwards (task 6). In this case, the "where" decision depends on the rooms already assigned to courses by the registrar. If a pandemic occurs before the registrar has assigned rooms to courses, our approach can also decide when and where each course meets in a socially distanced manner with total freedom to choose among rooms (task 4A). Even if a pandemic is nowhere in sight, our approach can create a backup schedule for an entire campus—performing task 4A or 6—with minimal administrative effort. Lastly, our approach can perform room assignment in normal times (task 4). More details about the generality of our approach are provided at the end of Section 4.

*1.2. Simultaneous attendance*

In this study, we propose a university course scheduling framework that enables all course sections—even the largest—to have a significant classroom learning component during a pandemic. Our approach does not allow rotating attendance and instead uses the idea of *simultaneous attendance* in which all students in a course section meet in multiple rooms at the same time, but less often than in a normal semester. At each *mass meeting*, students socially distance in each room. The instructor teaches in one room, and an image of the instructor is displayed in the other rooms simultaneously. Mass meetings provide opportunities for group work and in-person midterm exams that would not otherwise exist.

Figure 2 illustrates the idea of simultaneous attendance for a simple example with 8 course sections, 4 classrooms, and a 2-week semester. The top portion shows the schedule in normal times without social distancing. Four sections (A, B, C, D) meet from 10-11:30 on Monday and Wednesday, and four sections (E, F, G, H) meet from 11:30-1:00 on Monday and Wednesday. The classroom where each section meets is also shown. The bottom of the figure shows a schedule with simultaneous attendance that might be used during a pandemic if classroom capacities shrink by a factor of four. Here, each section occupies four times as many rooms, one-fourth as often. For example, Section A has one mass meeting in rooms 1-4 from 10-11:30 on Monday in week 1 and three online meetings at the other times when it normally meets. Other sections are analogous. In the pandemic schedule, all students attend each meeting for every section either in a classroom (as shown) or online (not shown).

Simultaneous attendance offers several benefits including opportunities for in-person midterm exams, but it requires a more complex scheduling approach than rotating attendance. With rotating attendance, no decision is made regarding when a course section occupies a room; it occupies a room as in a normal semester. Usually, no decision is made about where the section meets; it generally remains in the room initially assigned to it by the registrar. For example, the schedule at the top of Figure 2 can still be used during a pandemic if every course section is taught in hybrid mode with rotating attendance.

With simultaneous attendance during a pandemic, the scheduling is more complex. Each course

**Figure 2.** Schedule for 8 course sections, 4 rooms, and 2 weeks during normal times (a) and during a pandemic (b)

(a) Schedule for 8 course sections (A-H) during normal times without social distancing.

| | WEEK 1 | | | | WEEK 2 | | | |
|---|---|---|---|---|---|---|---|---|
| | Mon | | Wed | | Mon | | Wed | |
| Room | 10:00 | 11:30 | 10:00 | 11:30 | 10:00 | 11:30 | 10:00 | 11:30 |
| 1 | **A** | **E** | **A** | **E** | **A** | **E** | **A** | **E** |
| 2 | **B** | **F** | **B** | **F** | **B** | **F** | **B** | **F** |
| 3 | **C** | **G** | **C** | **G** | **C** | **G** | **C** | **G** |
| 4 | **D** | **H** | **D** | **H** | **D** | **H** | **D** | **H** |

(b) Classroom schedule with simultaneous attendance and social distancing during a pandemic. Classroom capacities are down 75%, so each section occupies 4 rooms, one-fourth as often. Three online meetings for each section are not shown.

| | WEEK 1 | | | | WEEK 2 | | | |
|---|---|---|---|---|---|---|---|---|
| | Mon | | Wed | | Mon | | Wed | |
| Room | 10:00 | 11:30 | 10:00 | 11:30 | 10:00 | 11:30 | 10:00 | 11:30 |
| 1 | **A** | **E** | **B** | **F** | **C** | **G** | **D** | **H** |
| 2 | **A** | **E** | **B** | **F** | **C** | **G** | **D** | **H** |
| 3 | **A** | **E** | **B** | **F** | **C** | **G** | **D** | **H** |
| 4 | **A** | **E** | **B** | **F** | **C** | **G** | **D** | **H** |

section usually occupies more than one room when it meets, so it is unlikely that all its meetings can be in person; some will be moved online. So, there is a *when* and *where* decision for each course section. Regarding when, we require the mass meetings to take place during a *subset* of the section's normal, in-person meeting times. Regarding where, we require the mass meetings to occupy a *fixed superset* of rooms, i.e., a set of rooms that does not change from one mass meeting to the next and includes the room assigned to the section by the registrar (if there is one). For example, at the top of Figure 2, Section A meets (in person) MW from 10-11:30 in room 1, so we require all its mass meetings during a pandemic to include room 1 and to take place from 10-11:30 on Monday or Wednesday. The schedule at the bottom of Figure 2 fulfills these "subset/superset" requirements for each course section.

*1.3. High-precision course scheduling*

Whereas most course scheduling methods create a repeating, weekly schedule, we take a different approach: we create a *high-precision*, nonrepeating schedule that explicitly considers each day of the semester. We do so for three reasons. First, some course sections have only a few in-person meetings a semester. In the case study described in Section 7, there are 52 *irregular* course sections, most of which meet in person less than once per week during a normal semester. The only way to schedule in-person meetings for these sections without wasting large amounts of space is to abandon a repeating weekly schedule. Second, many course sections meet once per week. If we require simultaneous attendance during a pandemic, many of these sections will meet in person less than once per week. The only way to schedule mass meetings for these sections without wasting large amounts of space is to abandon a repeating weekly schedule. Third, even if all course sections meet multiple times per week, the diversity of room capacities and section enrollments, meeting frequencies, and meeting times—many of which partially overlap—make it difficult to identify a regular, weekly time when every section's mass meetings take place during a pandemic. One of the goals of this work is to allow every course section to have a meaningful amount of in-person instruction during a pandemic, regardless of the number of sections or the diversity in section enrollments, room capacities, and meeting times/frequencies. These considerations motivate our use of a nonrepeating schedule that explicitly considers each day of the semester.

*1.4. Contributions*

In this paper, we present a novel approach to university course scheduling that enables all course sections—even the largest—to have a significant classroom learning component during a pandemic. Our approach is guided by the principles of fairness and simultaneous attendance which have received little attention in the pandemic-driven course scheduling literature. Fairness is achieved by ensuring that a basic fraction (e.g., 25%) of the pre-pandemic, in-person instruction in every course section remains in a classroom during a pandemic. Unlike previous studies, we do not allow rotating attendance and instead require simultaneous attendance in which all students in a section meet in 1-5 rooms at the same time, but less often than in a normal semester. Simultaneous attendance provides opportunities for group work and in-person midterm exams that would not otherwise exist. Instead of a weekly timetable, we create a high-precision schedule that considers each day of the semester. This allows us to handle all exceptions at [*StateXYZ*]'s second largest university including pre-scheduled, irregular midterm exams; sections with irregular in-person meetings; and temporary room closures for medical/vaccination activities.

We call the problem we address the *high-precision university course scheduling problem* (HP-UCSP). The word "pandemic" is not in the name because our approach also handles room assignment in normal times. We model this problem as an integer program and develop a fast heuristic algorithm that generates solutions to large problem instances in an hour.

Our ability to quickly tackle a large, real problem shows that our approach is *practical* and *scalable*. It is *flexible* by considering seven optimality criteria that are weighted according to user preferences. Our approach is also *general* as it performs three different tasks shown in Figure 1: regular room assignment (task 4 in Fig. 1); course scheduling during a pandemic without room restrictions (task 4A); and course scheduling during a pandemic with room restrictions (task 6). Finally, our approach is *friendly* to administrators because it places minimal demands on their time. High-precision schedules that prepare a campus for various pandemic possibilities can be created with minimal administrative input.

Our final contribution is to develop an expanded taxonomy of course delivery modes to facilitate future research. This taxonomy includes the modes considered by previous researchers, new ones introduced here, and several additional modes not yet explored.

This paper is organized as follows. Section 2 discusses the literature. In Section 3, we introduce an expanded taxonomy of course delivery modes. The high-precision university course scheduling problem (HP-UCSP) is formally described in Section 4. A math model of the HP-UCSP is presented in Section 5. In Section 6, we describe a heuristic algorithm for the HP-UCSP. In Section 7, we present the results of a case study of the [*UniversityXYZ*] campus in both pandemic and normal times. Concluding remarks are made in Section 8.

## 2. Literature review

Education has been a focal point of operations research for quite some time. According to Johnes (2015), operations research has been used in an educational context since the 1960s. Areas of application include planning and resource allocation for the educational sector, efficiency and performance measurement, routing, and scheduling. Of particular concern here is the literature on scheduling.

Research on educational scheduling includes two areas: school timetabling—which concerns primary and secondary schools—and university course scheduling. Pillay (2014) reviews the work on the school timetabling problem (STP). Teacher schedules and workloads play a central role in this problem.

University course scheduling generally considers more courses, rooms, students, instructors, and timeslots than school timetabling. The instructor of each course is fixed a-priori, and the main goal is to maximize student satisfaction. Most studies on university course scheduling address one of six problems: curriculum-based course timetabling (CB-CTT), post-enrollment course timetabling (PE-CTT), examination timetabling (ETT), strategic scheduling, timetabling by itself, and room assignment by itself.

Much of the research on CB-CTT, PE-CTT, and ETT originated with the Second International Timetabling Competition in August 2007 and its predecessor in 2002 (McCollum et al., 2010). An overview of CB-CTT research is provided by Bettinelli et al. (2015). In CB-CTT, tasks 3-4 in Figure 1 are done simultaneously: the university decides a weekly assignment of courses to rooms and time periods before students enroll in courses. The timetable must ensure that courses are only taught when instructors are available and courses belonging to the same curriculum (e.g., year 3 physics) or taught by the same instructor are not scheduled at the same time. There are also several soft constraints.

In PE-CTT, the university decides a weekly assignment of courses to rooms and time periods after students enroll in courses. Conflicts in PE-CTT are determined by students who individually enroll in particular courses, not by the curricula published by the university. In this approach, tasks 1-5 in Figure 1 are performed in the order 1-2-5-3-4 with tasks 3-4 performed simultaneously. Rudová et al. (2011) describe a method for performing PE-CTT at Purdue University that considers the large lecture courses, small lecture courses, and laboratory courses in sequential fashion.

Examination timetabling (ETT) resembles PE-CTT but allows exams to share rooms. Each exam needs to be assigned to a single room and period, and there are a variety of exam-specific considerations. For example, giving one exam later than another may be a requirement (McCollum et al., 2012).

Strategic scheduling is a broad area that includes relatively few papers. Some studies in this area consider how many classrooms and time periods—and which rooms and periods—are necessary (Lindahl et al., 2018). Other studies consider term planning (Khamechian and Petering, 2022).

Another approach to course timetabling and room assignment is to deal with them separately. This is a less student-focused approach than CB-CTT and PE-CTT as it does not consider distances that students travel between lectures. This approach may be used at large universities where courses are plentiful and elective courses form a major part of the curriculum. For example, Phillips et al. (2015) consider a large, practical room assignment problem at the University of Auckland.

Research on course scheduling during a pandemic began only recently. Barnhart et al. (2022) present math models and solution approaches for two problems faced by the Massachusetts Institute of Technology (MIT) as it prepared for the fall 2020 semester during the early stages of the COVID-19 pandemic. The first is a strategic term planning problem that helped decide (a) whether MIT should move from a two-semester to a three-semester calendar and (b) which courses should be offered in each term in each case. The second is a pandemic-driven, combined course timetabling, room assignment, and course modality selection problem with the goal of allowing students to attend lectures in person as much as possible despite reduced classroom capacities. In addressing the second problem, the authors introduce the ideas of teaching in multiple rooms and rotating attendance. They produce a weekly schedule using hierarchical optimization that was used in fall 2020 by the Sloan School of Management that allowed 68% of Sloan students to have in-person learning in at least half of their courses.

Johnson and Wilson (2022) consider a combined room assignment and course modality selection problem faced in Spears School of Business at Oklahoma State University prior to fall 2020. Course times are not changed, and the main goal is to maximize the number of courses delivered in person or in hybrid mode (with social distancing). The authors also introduce the idea of rotating attendance and use hierarchical optimization and an Excel-based solver to make a weekly schedule that was used in fall 2020.

Navabi-Shirazi et al. (2022) develop a math model and solution approach for a combined room assignment and course modality selection problem faced by Georgia Tech (GT) as it prepared for fall 2020. They consider a large problem with 2249 course sections and introduce a *hybrid touch point* delivery mode with rotating attendance in which each student attends in person less than once a week. Objectives include maximizing the satisfaction of modality preferences and previous room assignments and maximizing student-hours in the classroom (with social distancing). Hierarchical optimization is used to create a campus-wide schedule for GT in which 87% of course mode preferences are satisfied and 15.5% more student-hours are in a classroom than if no rooms are reassigned. Additional analyses show how (a) centralized versus decentralized planning, (b) the COVID-19 room capacity reduction percentage, and (b) the number of touch points needed for hybrid touch point mode affect the results.

Gore et al. (2022) develop two math models—a three-cohort model (3CM) and once-a-week model (OWM)—used by Clemson University to prepare for fall 2020. Assuming all courses are hybrid with rotating attendance, the models minimize the (a) number of unique students the average student shares a classroom with during the semester and (b) number of days an average student is on campus each week. The 3CM divides the entire student body into three cohorts which specify the three cohorts used for rotating attendance in each course. The OWM assumes each student attends a course once per week in person and the number of cohorts for a course equals the number of times it meets per week. The OWM decides which students are in each cohort in each course, and it was implemented in fall 2020.

Moug et al. (2022) develop two math models to help universities prepare for, and respond to, a pandemic. The first uses a conflict matrix to determine a classroom's capacity with social distancing requirements. The second model decides new room assignments and course modalities in the event of a mid-semester pandemic. It is solved via hierarchical optimization and used to create a weekly schedule of 1190 engineering courses at the University of Michigan in fall 2021. Results show that, under a 3-foot social distancing requirement, the model increases the fraction of lectures that students attend in person from 62% to 85% and the fraction of courses delivered in one room in person from 23% to 76%.

In this paper we propose an approach for scheduling university courses during a pandemic that differs from the preceding research in several ways. First, we do not allow rotating attendance and instead require all classroom attendance to consist of mass meetings in which all students in a course section meet in 1-5 rooms at the same in a socially distanced manner. Second, we ensure that a certain fraction of the instruction in every course section occurs in the classroom. Third, we schedule across all days of the semester instead of a single, repeating week. This allows us to accommodate pre-scheduled midterm exams and course sections with irregular in-person meetings. Fourth, we do not use hierarchical optimization and instead use a weighted objective function. Our approach handles all exceptions at [*StateXYZ*]'s second largest university and allows at least 25% of the instruction in every section, and more than 49% of all instruction across the entire campus, to be in a classroom.

## 3. Expanded taxonomy of course delivery modes

In this section, we propose an expanded taxonomy of course delivery modes which includes the concepts proposed thus far and new ideas for future research. The taxonomy is shown in Table 1.

Course delivery refers to how students receive oral instruction from the instructor of a particular course section. A course section (i.e., *section*) is either part of a course or equivalent to a course. For example, a large course with 2 lecture sections, 17 discussion sections, and 15 lab sections translates to 34 sections. A small course with 1 lecture section and no discussion or lab sections translates to 1 section.

Course delivery can be broadly categorized into three *modes*—in person (IP), online (OL), and hybrid (HY). We propose that these modes be defined as follows. In mode IP, all students attend all lectures for the course section—which we call *meetings*—in a classroom. In mode OL, all students attend all meetings online. In mode HY, some attendance during the semester is in person and some is online.

Note that our definition of mode IP does not state that all students attend all lectures in the same room. Indeed, it is possible for students to be split across multiple rooms (Barnhart et al., 2022). Hence, we propose that the delivery mode for an in-person section be specified using the abbreviation IP-*r* where *r* indicates how many rooms are used.

Mode HY can be delivered in at least three *formats*—with fixed (-F), rotating (-R), or simultaneous (-S) attendance. These formats are indicated using the abbreviations HY-F, HY-R, and HY-S. In mode HY-F, each meeting has a classroom audience and a remote audience, and the method of attendance is fixed for each student. Some students attend all meetings in person; all other students attend all meetings remotely. As Navabi-Shirazi et al. (2022) have shown, a hybrid section may have rotating

**Table 1.** Taxonomy of course delivery modes

| Delivery Mode Abbreviation | Description |
|---|---|
| IP-*r* | In-person. All students attend all meetings in a classroom, and *r* classrooms are used. |
| OL | Online. All students attend all meetings remotely. No classrooms are used. |
| HY-F-*r* | Hybrid with fixed attendance and *r* classrooms. Each meeting has a classroom and remote audience. The method of attendance—in person or remote—is fixed for each student throughout the semester. |
| HY-R-W-*r* | Hybrid with a rotating, weekly attendance roster in *r* classrooms. Each meeting has a classroom and remote audience. Each student attends at least 1 meeting per week in person and attends other meetings remotely. Attendance roster repeats weekly. |
| HY-R-N-*r* | Hybrid with a rotating, nonweekly attendance roster in *r* classrooms. Each meeting has a classroom and remote audience. Attendance roster repeats regularly but not weekly. |
| HY-S-W-*r* | Hybrid with simultaneous attendance in *r* classrooms on a repeating, weekly basis. Each meeting has a classroom audience or remote audience but not both. Each student attends at least 1 meeting per week in a classroom and attends other meetings remotely. Schedule repeats weekly. |
| HY-S-N-*r* | Hybrid with simultaneous attendance in *r* classrooms on a repeating, nonweekly basis. Each meeting has a classroom or remote audience but not both. Schedule repeats regularly but not weekly. |
| HY-S-I-*r* | Hybrid with simultaneous attendance in *r* classrooms on an irregular subset of days. Each meeting has a classroom or remote audience but not both. All students attend in person on specific days identified in advance and attend other meetings remotely. Schedule does not repeat regularly. |

(-R) attendance on a weekly (-W) or a nonweekly (-N) basis. In the former case, which we call HY-R-W, the meeting *attendance roster*—the list of students who attend each meeting in a classroom—repeats weekly. The other case, HY-R-N, includes two possibilities: students attending in person less than once per week on average (e.g., once every two weeks) or more than once per week on average (e.g., three times every two weeks). For the sake of brevity, we use HY-R-N for both.

In this work, we introduce the hybrid delivery mode with simultaneous attendance (HY-S). Here, each meeting is attended by a classroom audience or online audience but not both. When the meeting is in the classroom, we call it a *mass meeting*. If the mass meetings repeat weekly, the abbreviation HY-S-W is used. If they happen at a regular frequency but not weekly (e.g., once every two weeks), the term HY-S-N is used. If the meetings happen on an irregular set of days without obvious repetition, the term HY-S-I applies. For example, consider a course section that normally meets MWF from 9-10. If the section's mass meetings during a pandemic happen every Monday and Friday from 9-10, the term HY-S-W applies. If they happen every third Wednesday from 9-10, the term HY-S-N applies. If they happen from 9-10 on (M, F, W, W, M, W) in weeks (1, 1, 5, 7, 9, 12) of the semester, the term HY-S-I applies.

For HY sections, the number of classrooms used ($r$) can be indicated by a suffix *-r* at the end of the modality abbreviation. For example, HY-R-W-2 indicates the option utilized by Barnhart et al. (2022) in which a hybrid section with rotating weekly attendance is taught in two classrooms.

Using the above taxonomy, Table 2 shows the delivery modes examined in the recent literature on pandemic-driven course scheduling. Note that modes HY-S-W and HY-S-N have not been explored.

**Table 2.** Course delivery modes considered in the literature

| Delivery Mode | Barnhart et al. (2022) | Johnson and Wilson (2022) | Navabi-Shirazi et al. (2022) | Gore et al. (2022) | Moug et al. (2022) | This work |
|---|---|---|---|---|---|---|
| IP-1 | x | x | x |  | x | x |
| IP-2 | x |  |  |  | x | x |
| IP-$r$ ($r \geq 3$) |  |  |  |  |  | x |
| OL | x | x | x |  | x |  |
| HY-F |  |  |  |  |  |  |
| HY-R-W-1 | x | x | x | x | x |  |
| HY-R-W-2 | x |  |  |  |  |  |
| HY-R-W-$r$ ($r \geq 3$) |  |  |  |  |  |  |
| HY-R-N-1 |  | x | x | x | x |  |
| HY-R-N-2 |  |  |  |  |  |  |
| HY-R-N-$r$ ($r \geq 3$) |  |  |  |  |  |  |
| HY-S-W |  |  |  |  |  |  |
| HY-S-N |  |  |  |  |  |  |
| HY-S-I-1 |  |  |  |  |  | x |
| HY-S-I-2 |  |  |  |  |  | x |
| HY-S-I-$r$ ($r \geq 3$) |  |  |  |  |  | x |

## 4. Formal problem description

Our problem is called the high-precision university course scheduling problem (HP-UCSP). Consider a university that has completed tasks 1-5 of the course scheduling process for an upcoming semester (Fig. 1). There is no pandemic, and the university plans to offer $S$ course sections, including sections to be delivered in person and in hybrid mode but not online. Let $E_s$ be the number of students enrolled in, or expected to be enrolled in, section $s$.

The semester has $W$ weeks and $D$ academic days. If instruction occurs five (six) days a week, $D = 5W$ ($6W$). Let $N_{sw}$ be the number of in-person *meetings* for section $s$ that take place in week $w$ and $CN_{sw}$ be the number of such meetings that take place in weeks 1 to $w$ combined. These values account for holidays and hybrid sections. For example, if section $s$ meets in person three times per week, $N_{sw} = 3$ for all $w$ except possibly weeks with holidays and the first/last weeks of the semester. If hybrid section $s$ meets once a week but only every fourth week in person, $N_{sw} = 1$ every fourth week and 0 otherwise.

The academic week is divided into $T$ *timeslots* of equal duration, and the instruction of any section on any day fits into an integral number of consecutive timeslots. If all sections' meetings start at :00 or :30 after the hour, each timeslot lasts 30 minutes. If there is more variety in meeting start times, each timeslot is shorter. The value of $T$ depends on the timeslot duration, number of days of instruction each week, and earliest start and latest finish time of any meeting on any day. If each timeslot is 30 minutes and if instruction occurs 6 days per week during an 8:00-21:00 window, $T = 156$ (= 6*26).

Instruction occurs during $M$ weekly *meeting times*. Each meeting time is an uninterrupted period during the week when the instruction for one or more sections occurs. For example, three sections that meet MWF from 9-10, MW from 9-10, and MW from 9-10:30 use five unique meeting times: M 9-10, W 9-10, F 9-10, M 9-10:30, and W 9-10:30. Binary parameter $I_{sm}$ equals 1 if section $s$ has an in-person meeting during meeting time $m$ in any week and 0 otherwise. Binary parameter $O_{mt}$ equals 1 if meeting time $m$ includes timeslot $t$ and 0 otherwise.

Sections with irregular in-person meetings are modeled using binary parameter $K_{sd}$ which equals 1 (0) if section $s$ does not have (does have) an in-person meeting on day $d$. For example, if hybrid section $s$ meets twice in person during the semester, then $K_{sd} = 1$ for all $d$ except the two days of the in-person meetings. The days and times of all in-person meetings for section $s$ are given by $I_{sm}$ and $K_{sd}$ combined.

There are $R$ classrooms (i.e., *rooms*) on campus, and binary parameter $InitRoom_{rs}$ equals 1 if room $r$ has been assigned to section $s$ and 0 otherwise. If rooms have not been assigned, $InitRoom_{rs}$ equals 0 for all $r$ and $s$. If a room has been assigned to every section, $S$ of these parameters equal 1.

The HP-UCSP is particularly relevant if a pandemic occurs and classroom capacities suddenly decrease because social distancing is required in the classroom. Given this crisis, the university replans the semester as follows. Going forward, each section will have a limited number of *mass meetings* during the semester when all students in the section meet in one or more rooms at the same time, with all other meetings online. The times of the mass meetings will be a subset of the times when the section normally meets, and the rooms assigned to the section will include the room already assigned to it (if any).

Due to irregularities already built into the schedule, mass meetings will not occur in repeating fashion (e.g., weekly or biweekly). Instead, they will be scheduled across all $W$ weeks and $D$ days of the semester. The university will not allow any section to go fully online. In pursuit of fairness, it declares that at least $MinFraction_s$ of the in-person meetings that were planned for section $s$ prior to the pandemic will still be in person during the pandemic ($0 < MinFraction_s \leq 1$ for all $s$).

The replanning effort considers the following information. First, $C_r$ is the adjusted capacity of room $r$ with social distancing requirements. Binary parameter $J_{rs}$ equals 1 if the equipment in room $r$ is compatible with section $s$ and 0 otherwise. Some rooms may not be available on some days due to construction, flooding, or medical/vaccination activities, so binary parameter $NAvail_{rd}$ equals 1 (0) if room $r$ is not (is) available on day $d$.

The university's overall task is to create a new schedule for each section during the pandemic that specifies the room(s) where, and (meeting time, week) pairs when, its mass meetings take place. The rooms where each section's mass meetings occur must remain fixed for the entire semester.

The pandemic schedule for each section is judged according to seven criteria. (1) The number of rooms used for the section should be minimized. (2) The distance between the rooms used for the section should be minimized. (3) The rooms used for the section should be in buildings that are preferred by the *organization* (i.e., department) that teaches it. (4) The number of wasted (i.e., unutilized) seats in the section's room assignment should be minimized. (5) The number of mass meetings for the section should be maximized. (6) The mass meetings should be spread throughout the semester, not packed into just a few weeks' time. (7) The fraction of the in-person meetings that were planned for section $s$ prior to the pandemic and will still be in person during the pandemic should be at least $MinFraction_s$. Finally, because less experienced students need more in-person contact than more experienced students, there should be a preference for better schedules in lower-level (e.g., first year) courses than higher-level courses.

Note that the HP-UCSP applies to the creation of a schedule before, or during, a semester. It also applies to room assignment in normal times (task 4 in Figure 1) if exactly one room is assigned to each section; all $InitRoom_{rs} = 0$; $C_r$ is the normal (i.e., non-pandemic) capacity of room $r$; and $MinFraction_s = 1.0$ for all $s$. With this in mind, we now introduce an integer programming formulation of the problem.

## 5. Mathematical model

Our mathematical model, model HP-UCSP, is shown in Tables 3-6 and equations 1-37. Table 3 shows the indices, and Table 4 shows the input data used in the model. Table 5 shows the decision variables and their domains. Table 6 shows the model's adjustable parameters. These can be changed by the user to pursue different objectives. Equation 1 is the objective function, and equations 2-37 are the constraints. The items in Tables 3-6 are generally sequenced in order of their appearance in equations 1-37 and are explained as we discuss the objective function and constraints on the following pages.

The model is an integer program with three main decision variables shown at the top of Table 5. The first, $X_{smw}$, equals 1 if a mass meeting for section $s$ is scheduled during meeting time $m$ in week $w$ and 0 otherwise. The second, $Y_{rs}$, equals 1 if room $r$ is assigned to section $s$ and 0 otherwise. The third, $Z_{rsmw}$, combines the information in the $X$ and $Y$ variables; it equals 1 if room $r$ hosts a mass meeting for section $s$ during meeting time $m$ in week $w$ and 0 otherwise.

Twenty auxiliary decision variables, shown at the bottom of Table 5, help to enforce constraints and translate the values of the main decision variables into objective penalties. The first eight of these appear in the objective function; the last twelve do not.

The objective (1) is to minimize a sum of seven objective *components*, each weighted by a different factor $\alpha_1$ to $\alpha_7$. Components 1-4 (5-7) are impacted by where (when) each section meets. Each component is a sum of penalties for the $S$ sections, with the penalty for section $s$ weighted by a factor that depends on its importance (*Importance$_s$*) as determined by its course level, its enrollment ($E_s$), the number of in-person meetings it has during a normal semester ($CN_{sW}$) (note the capital $W$), and the duration of its meetings (*Duration$_s$*). The term $E_s*CN_{sW}*Duration_s$, the number of student-timeslots of in-person attendance in section $s$ during a normal semester, appears in five objective components.

Penalty 1 equals the number of rooms assigned to section $s$ (*NumRooms$_s$*) minus 1. Penalty 2 is *DistPenalty$_s$*, an aggregate measure of the distance between rooms assigned to section $s$. It equals 0 if

**Table 3.** Indices in model HP-UCSP

| Index | Description |
|---|---|
| s | Section (i.e., course section) (1 to $S$) |
| r, q | Room (i.e., classroom) (1 to $R$) |
| m | Meeting time; uninterrupted weekly period when the instruction for one or more sections occurs (1 to $M$) |
| w, k | Week (1 to $W$) |
| d | Day of instruction during the semester (1 to $D$) |
| t | Timeslot; period of fixed duration (e.g., 30 minutes) during the academic week (1 to $T$) |
| b, c | Building (1 to $B$) |
| f | Floor within a building (0 to $F$; 0 is the basement) |
| g | Organization that offers courses (1 to $G$) |
| l | Level (i.e., course level) (0 to $L$; 0 represents a remedial course) |

**Table 4.** Input data for model HP-UCSP (integer-valued unless otherwise noted)

| Data Item | Description |
|---|---|
| $S$ | Number of sections with an in-person learning component in a normal (non-pandemic) semester |
| $R$ | Number of rooms where instruction can occur |
| $M$ | Number of unique weekly meeting times when the $S$ sections are taught |
| $W$ | Number of weeks in the semester |
| $D$ | Number of academic days in the semester |
| $T$ | Number of timeslots of equal duration (e.g., 30 minutes) in each academic week |
| $B$ | Number of buildings where instruction can occur |
| $F$ | Highest floor number (in any building) where instruction can occur |
| $G$ | Number of organizations (i.e., departments) that offer courses |
| $L$ | Highest possible course level (e.g., 9 for a 900-level course) |
| $L_s$ | Course level of section $s$ (e.g., 1 for a 100-level course, 9 for a 900-level course) |
| $E_s$ | Enrollment in section $s$. (Number of students enrolled in, or expected to be enrolled in, section $s$.) |
| $N_{sw}$ | Number of in-person meetings for section $s$ that occur in week $w$ in a normal semester |
| $CN_{sw}$ | Number of in-person meetings for section $s$ that occur in weeks 1 to $w$ combined in a normal semester |
| $Duration_s$ | Number of timeslots used by each meeting of section $s$ |
| $I_{sm}$ | = 1 if section $s$ has an in-person meeting during meeting time $m$ (in any week) in a normal semester<br>= 0 otherwise |
| $Q_{mwd}$ | = 1 if meeting time $m$ in week $w$ falls on day $d$<br>= 0 otherwise |
| $K_{sd}$ | = 1 if section $s$ is not allowed to have an in-person meeting on day $d$<br>= 0 if section $s$ is allowed to have an in-person meeting on day $d$ |
| $H_d$ | = 1 if day $d$ is a holiday (when no instruction occurs)<br>= 0 otherwise |
| $C_r$ | Capacity of room $r$ (during a pandemic or otherwise) |
| $J_{rs}$ | = 1 if room $r$ can be used for section $s$<br>= 0 otherwise |
| $InitRoom_{rs}$ | = 1 if room $r$ is already assigned to section $s$ (e.g., by the university registrar)<br>= 0 otherwise |
| $O_{mt}$ | = 1 if meeting time $m$ includes timeslot $t$<br>= 0 otherwise |
| $NAvail_{rd}$ | = 1 if room $r$ is not available on day $d$ (e.g., due to construction, flooding, vaccination events, etc.)<br>= 0 if room $r$ is available on day $d$ |
| $IB_{rb}$ | = 1 if room $r$ is in building $b$<br>= 0 otherwise |
| $Dist_{bc}$ | Distance (in meters) between the centroids of buildings $b$ and $c$ (real, $\geq 0$) |
| $IF_{rf}$ | = 1 if room $r$ is on floor $f$ of its building<br>= 0 otherwise |
| $NAdj_{qr}$ | = 1 if rooms $q$ and $r$ are not adjacent, not on the same floor, or not in the same building<br>= 0 if rooms $q$ and $r$ are adjacent (and on the same floor of the same building) |
| $G_s$ | Organization that teaches section $s$ |
| $B_r$ | Building where room $r$ resides |
| $P_{gb}$ | Building preference penalty when organization $g$ teaches in building $b$. (The most preferred $g$-$b$ combinations have a value of 0. Less preferred $g$-$b$ combinations have higher values.) |
| $FirstW_s$ | First week when section $s$ has an in-person meeting in a normal semester |
| $LastW_s$ | Last week when section $s$ has an in-person meeting in a normal semester |

**Table 5.** Decision variables in model HP-UCSP

| Decision Variable | Description |
|---|---|
| $X_{smw}$ | = 1 if a mass meeting for section $s$ is scheduled during meeting time $m$ in week $w$ <br> = 0 otherwise (binary) |
| $Y_{rs}$ | = 1 if room $r$ is assigned to section $s$ <br> = 0 otherwise (binary) |
| $Z_{rsmw}$ | = 1 if room $r$ hosts a mass meeting for section $s$ during meeting time $m$ in week $w$ <br> = 0 otherwise (real, $\geq 0$) |
| $NumRooms_s$ | Number of rooms assigned to section $s$ (real, $\geq 1$) |
| $DistPenalty_s$ | Penalty for the distance between rooms assigned to section $s$ (real, $\geq 0$) |
| $PrefPenalty_s$ | Penalty for teaching section $s$ in nonpreferred buildings (real, $\geq 0$) |
| $WastedSeats_s$ | Number of empty seats in section $s$'s room assignment (real, $\geq 0$) |
| $NP_{sw}$ | Number of mass meetings for section $s$ scheduled in week $w$ during the pandemic (real, $\geq 0$) |
| $CNP_{sw}$ | Number of mass meetings for section $s$ scheduled in weeks 1 to $w$ combined during the pandemic (real, $\geq 0$) |
| $TimingPenalty_s$ | Penalty for the irregularity in the timing of section $s$'s mass meetings (real, $\geq 0$) |
| $U_s$ | = 1 if at least $MinFraction_s$ of section $s$'s meetings during the pandemic are in person <br> = 0 otherwise (binary) |
| $MaxBuildingDist_s$ | Maximum distance (in meters) between buildings where section $s$ is taught (real, $\geq 0$) |
| $NumBldgs_s$ | Number of buildings in which section $s$ is taught (real, $\geq 1$) |
| $FloorDist_s$ | Total distance between floors where section $s$ is taught, summed over all buildings (real, $\geq 0$) |
| $NumExtraFloors_s$ | Number of extra floors, beyond one in each building, used for teaching section $s$ (real, $\geq 0$) |
| $RoomNonadj_s$ | Number of nonadjacent room pairs in section $s$'s room assignment (real, $\geq 0$) |
| $BU_{bs}$ | = 1 if building $b$ is used for teaching section $s$ <br> = 0 otherwise (binary) |
| $FU_{fbs}$ | = 1 if floor $f$ in building $b$ is used for teaching section $s$ <br> = 0 otherwise (binary) |
| $HighFloor_{bs}$ | Highest floor in building $b$ where section $s$ is taught (real, $\geq 0$) |
| $LowFloor_{bs}$ | Lowest floor in building $b$ where section $s$ is taught (real, $\geq 0$) |
| $NumFloors_s$ | Number of unique floors (in one or more buildings) in which section $s$ is taught (real, $\geq 1$) |
| $YY_{qrs}$ | = 1 if rooms $q$ and $r$ are both assigned to section $s$ <br> = 0 otherwise (real, $\geq 0$) |
| $Diff_{sw}$ | Difference between the actual and prorated number of mass meetings for section $s$ that are scheduled in weeks 1 to $w$ combined (real, $\geq 0$) |

**Table 6.** Adjustable parameters in model HP-UCSP (integer-valued unless otherwise noted)

| Parameter | Description |
|---|---|
| $\alpha_1, \alpha_2, \ldots, \alpha_7$ | Weights for the seven objective components (real, $\geq 0$) |
| $\alpha_{21}, \alpha_{22}, \ldots, \alpha_{25}$ | Weights for the five subcomponents of objective component 2 (real, $\geq 0$) |
| $Exp$ | Exponent used in objective component 5 (real, $\geq 1$) |
| $ImportanceOfLevel_l$ | Importance of level $l$ sections (e.g., 5 = most important; 1 = least important) (real, $\geq 0$) |
| $Importance_s$ | Importance of section $s$ (based on the level of section $s$ ($L_s$) and the value above) |
| $MaxRooms$ | Maximum number of rooms that can be assigned to any section |
| $MaxDistPenalty_s$ | Maximum allowed "room distance penalty" for section $s$ (real, $\geq 0$) |
| $MaxPrefPenalty_s$ | Maximum allowed building preference penalty for section $s$ |
| $MaxWastedSeats_s$ | Maximum allowed number of wasted (i.e., empty) seats in section $s$'s room assignment |
| $MinFraction_s$ | Minimum fraction of section $s$'s meetings that must be in-person mass meetings (real, $\geq 0$) |
| $BigM$ | Large positive number |

only one room is assigned or if several rooms in very close proximity on the same floor of the same building are assigned. More details are provided in the discussion of constraints 13-26. Penalty 3 is *PrefPenalty$_s$*, the penalty for teaching section *s* in nonpreferred buildings. Penalty 4 is *WastedSeats$_s$*, the number of empty seats in section *s*'s room assignment. For example, there are 8 wasted seats if a section with enrollment 42 is assigned to two rooms each with capacity 25. Penalty 5 is the fraction of section *s*'s meetings that are online $((CN_{sW} - CNP_{sW})/CN_{sW})$, raised to the power *Exp* (note the capital *W*). For example, if a section has $CN_{sW} = 20$ in-person meetings during a normal semester but only $CNP_{sW} = 8$ (in-person) mass meetings during a pandemic, then 60% $(= (20 - 8)/20)$ of its instruction has moved online. This is then raised to the *Exp* power. A value of *Exp* > 1 encourages this percentage to be similar across different sections (due to the convexity of $y = x^{Exp}$), and a value of 1 does not. Penalty 6 is *TimingPenalty$_s$*, a measure of the irregularity in the timing of section *s*'s mass meetings. More details are provided in the discussion of constraints 34-36. Penalty 7 is $(1 - U_s)$, which equals 1 if less than *MinFraction$_s$* of section *s*'s meetings are in person and 0 otherwise. A high value of $\alpha_7$ is needed to prioritize this penalty.

Minimize:

$(\alpha_1) * \sum_{s=1}^{S}[Importance_s * E_s * CN_{sW} * Duration_s * (NumRooms_s - 1)] +$

$(\alpha_2) * \sum_{s=1}^{S}[Importance_s * E_s * CN_{sW} * Duration_s * DistPenalty_s] +$

$(\alpha_3) * \sum_{s=1}^{S}[Importance_s * E_s * CN_{sW} * Duration_s * PrefPenalty_s] +$

$(\alpha_4) * \sum_{s=1}^{S}[CN_{sW} * Duration_s * WastedSeats_s] +$

$(\alpha_5) * \sum_{s=1}^{S}[Importance_s * E_s * CN_{sW} * Duration_s * ((CN_{sW} - CNP_{sW})/CN_{sW})^{Exp}] +$

$(\alpha_6) * \sum_{s=1}^{S}[Importance_s * E_s * Duration_s * TimingPenalty_s] +$

$(\alpha_7) * \sum_{s=1}^{S}[Importance_s * E_s * CN_{sW} * Duration_s * (1 - U_s)]$

(1)

We discuss the constraints in groups. Constraints 2a, 2b, and 2c ensure the consistency of $X_{smw}$, $Y_{rs}$, and $Z_{rsmw}$. Constraints 2a and 2b ensure that $Z_{rsmw}$ is no greater than 1 and is 0 if $X_{smw}$ or $Y_{rs}$ are 0. On the other hand, (2c) requires $Z_{rsmw}$ to be 1 if $X_{smw}$ and $Y_{rs}$ are both 1. These three constraints allow $Z_{rsmw}$ to be a real variable whose value is either 0 or 1 depending on the values of the binary variables $X_{smw}$ and $Y_{rs}$.

| | | |
|---|---|---|
| $Z_{rsmw} \leq X_{smw}$ | for all *r, s, m*, and *w* | (2a) |
| $Z_{rsmw} \leq Y_{rs}$ | for all *r, s, m*, and *w* | (2b) |
| $Z_{rsmw} \geq X_{smw} + Y_{rs} - 1$ | for all *r, s, m*, and *w* | (2c) |

Constraints 3-5 impose restrictions on $X_{smw}$. Constraint 3 ensures that a section's mass meetings during a pandemic occur during its normal meeting times. Constraint 4 handles sections that have irregular in-person meetings in normal times, and (5) ensures there is no instruction on holidays.

| | | |
|---|---|---|
| $X_{smw} \leq I_{sm}$ | for all $s$, $m$, and $w$ | (3) |
| $X_{smw} * Q_{mwd} * K_{sd} = 0$ | for all $s$, $m$, $w$, and $d$ | (4) |
| $X_{smw} * Q_{mwd} * H_d = 0$ | for all $s$, $m$, $w$, and $d$ | (5) |

Constraints 6-10 impose restrictions on $Y_{rs}$. Constraint 6 ensures that the rooms assigned to a section have sufficient capacity to accommodate it. Constraint 7 ensures that only suitable rooms are assigned to each section. Constraint 8 requires room $r$ to be assigned to section $s$ if the registrar already assigned it to the section. Constraints 9-10 compute the number of rooms assigned to a section and ensure that it does not exceed *MaxRooms*.

| | | |
|---|---|---|
| $\sum_{r=1}^{R} C_r * Y_{rs} \geq E_s$ | for all $s$ | (6) |
| $Y_{rs} \leq J_{rs}$ | for all $r$ and $s$ | (7) |
| $Y_{rs} \geq InitRoom_{rs}$ | for all $r$ and $s$ | (8) |
| $NumRooms_s = \sum_{r=1}^{R} Y_{rs}$ | for all $s$ | (9) |
| $NumRooms_s \leq MaxRooms$ | for all $s$ | (10) |

Constraints 11-12 impose restrictions on $Z_{rsmw}$. Constraint 11 ensures that at most one section uses a room at any time, and (12) ensures that no instruction occurs in room $r$ on day $d$ if $NAvail_{rd} = 1$.

| | | |
|---|---|---|
| $\sum_{s=1}^{S} \sum_{m=1}^{M} Z_{rsmw} * O_{mt} \leq 1$ | for all $r$, $t$, and $w$ | (11) |
| $Z_{rsmw} * Q_{mwd} * NAvail_{rd} = 0$ | for all $r$, $s$, $m$, $w$, and $d$ | (12) |

Constraints 13-26 quantify the distance between rooms assigned to a section. Constraint 13 shows that *DistPenalty$_s$*, the second objective penalty, is a weighted sum of five terms. The first, *MaxBldgDist$_s$*, is the maximum distance in meters between buildings where section $s$ is taught. This is zero if section $s$ is taught in one building. The second, (*NumBldgs$_s$* – 1), is the number of extra buildings used for section $s$ beyond the first. The third, *FloorDist$_s$*, is the total distance between floors where section $s$ is taught, summed over each building where it is taught. For example, if section $s$ uses four rooms that are located on floors 1, 2, and 6 of building A and floor 1 of building B, *FloorDist$_s$* equals 5 (= 6 – 1). The fourth term, *NumExtraFloors*$_s$, is the number of extra floors, beyond one per building, used for section *s*. In the previous example, this is 2 for the two extra floors used in building A. The fifth term, *RoomNonadj*$_s$, is the number of nonadjacent room pairs in section *s*'s room assignment. Two rooms are considered adjacent if they are on the same floor and either (i) share a wall and have doors in the same hallway or (ii) have doors directly across from each other. For example, if a section is taught in four rooms on four different floors, this is 6 = (4 nCr 2). If it is taught in four rooms in a row along the same side of a hallway, this is 3. If it is taught in four rooms whose doors are clustered at one location in a hallway, this is 0.

Constraint 14 ensures that *DistPenalty*$_s$ does not exceed a given limit (possibly ∞). Constraint 15 uses *Y*$_{rs}$ to compute *BU*$_{bs}$, which indicates if building *b* is used for section *s*. Constraint 16 uses *BU*$_{bs}$ and parameter *Dist*$_{bc}$ to compute *MaxBldgDist*$_s$. Constraint 17 uses *BU*$_{bs}$ to compute *NumBldgs*$_s$. Constraint 18 uses *Y*$_{rs}$ to compute *FU*$_{fbs}$, which indicates if floor *f* in building *b* is used for section *s*. Constraints 19-22 use *FU*$_{fbs}$ to compute *HighFloor*$_{bs}$ and *LowFloor*$_{bs}$, the highest and lowest floors in building *b* where section *s* is taught, which are then used to compute *FloorDist*$_s$. Constraint 23 uses *FU*$_{bs}$ to compute *NumFloors*$_s$, and (24) computes *NumExtraFloors*$_s$. The values of *BU*$_{bs}$, *MaxBldgDist*$_s$, *NumBldgs*$_s$, *FU*$_{fbs}$, *FloorDist*$_s$, *NumFloors*$_s$, and *NumExtraFloors*$_s$ are not accurate at all feasible solutions, but they are correct at every optimal solution due to their contribution to an objective function that is minimized. Constraints 25a, 25b, and 25c ensure the consistency of variables *YY*$_{qrs}$, *Y*$_{qs}$, and *Y*$_{rs}$. They allow *YY*$_{qrs}$, which indicates if rooms *q* and *r* are both used for section *s*, to be real instead of binary. Constraint 26 uses *YY*$_{qrs}$ and parameter *NAdj*$_{qr}$ to compute *RoomNonadj*$_s$. A detailed inspection of building floor plans is needed to determine the value of *NAdj*$_{qr}$ which equals 1 (0) if rooms *q* and *r* are not (are) adjacent.

$$DistPenalty_s = (\alpha_{21})(MaxBldgDist_s) + \\ (\alpha_{22})(NumBldgs_s - 1) + \\ (\alpha_{23})(FloorDist_s) + \\ (\alpha_{24})(NumExtraFloors_s) + \\ (\alpha_{25})(RoomNonadj_s) \quad \text{for all } s \quad (13)$$

$$DistPenalty_s \leq MaxDistPenalty_s \quad \text{for all } s \quad (14)$$

$$BU_{bs} * MaxRooms \geq \sum_{r=1}^{R} IB_{rb} * Y_{rs} \quad \text{for all } b \text{ and } s \quad (15)$$

$$MaxBldgDist_s \geq Dist_{bc} - (BigM)*(2 - BU_{bs} - BU_{cs}) \quad \text{for all } (b, c) \text{ s.t. } b < c \text{ and all } s \quad (16)$$

$$NumBldgs_s = \sum_{b=1}^{B} BU_{bs} \quad \text{for all } s \quad (17)$$

$$FU_{fbs} * MaxRooms \geq \sum_{r=1}^{R} IF_{rf} * IB_{rb} * Y_{rs} \quad \text{for all } f, b, \text{ and } s \quad (18)$$

$$HighFloor_{bs} \geq f - (F)*(1 - FU_{fbs}) \quad \text{for all } f, b, \text{ and } s \quad (19)$$

$$LowFloor_{bs} \leq f + (F)*(1 - FU_{fbs}) \qquad \text{for all } f, b, \text{ and } s \qquad (20)$$

$$HighFloor_{bs} \geq LowFloor_{bs} \qquad \text{for all } b \text{ and } s \qquad (21)$$

$$FloorDist_s = \sum_{b=1}^{B}(HighFloor_{bs} - LowFloor_{bs}) \qquad \text{for all } s \qquad (22)$$

$$NumFloors_s = \sum_{f=0}^{F}\sum_{b=1}^{B} FU_{fbs} \qquad \text{for all } s \qquad (23)$$

$$NumExtraFloors_s = NumFloors_s - NumBldgs_s \qquad \text{for all } s \qquad (24)$$

$$YY_{qrs} \leq Y_{qs} \qquad \text{for all } q, r, \text{ and } s \qquad (25a)$$

$$YY_{qrs} \leq Y_{rs} \qquad \text{for all } q, r, \text{ and } s \qquad (25b)$$

$$YY_{qrs} \geq Y_{qs} + Y_{rs} - 1 \qquad \text{for all } q, r, \text{ and } s \qquad (25c)$$

$$RoomNonadj_s = \sum_{q=1}^{R}\sum_{r=q+1}^{R} NAdj_{qr} * YY_{qrs} \qquad \text{for all } s \qquad (26)$$

Constraints 27-28 compute *PrefPenalty$_s$*, the third penalty term, and ensure it does not exceed a given limit (possibly ∞). On the right side of (27), $G_s$ is the organization that teaches section $s$, $B_r$ is room $r$'s building, and $P_{gb}$ is organization $g$'s level of undesirability for teaching in building $b$. This is higher for buildings that are farther from organization $g$'s office. Constraint 27 ensures that, at an optimal solution, *PrefPenalty$_s$* equals the building preference penalty of the worst building used for section $s$.

$$PrefPenalty_s \geq P_{G_s,B_r} - (BigM)*(1 - Y_{rs}) \qquad \text{for all } r \text{ and } s \qquad (27)$$

$$PrefPenalty_s \leq MaxPrefPenalty_s \qquad \text{for all } s \qquad (28)$$

Constraints 29-30 compute *WastedSeats$_s$*, the fourth penalty term, and ensure it does not exceed a user-defined limit (possibly ∞). Wasted seats are a system-wide loss, not a loss for the section that wastes them, so the factors *Importance$_s$* and $E_s$ are omitted from objective component 4 (see Eq. 1).

$$WastedSeats_s = (\sum_{r=1}^{R} C_r * Y_{rs}) - E_s \qquad \text{for all } s \qquad (29)$$

$$WastedSeats_s \leq MaxWastedSeats_s \qquad \text{for all } s \qquad (30)$$

Constraints 31-32 compute *CNP$_{sw}$*, the number of mass meetings for section $s$ in weeks 1 to $w$ combined, which is used in penalty 5. Constraint 33, which is optional, requires at least *MinFraction$_s$* of section $s$'s meetings to be in person (note the capital $W$).

$NP_{sw} = \sum_{m=1}^{M} X_{smw}$ for all $s$ and $w$ (31)

$CNP_{sw} = \sum_{k=1}^{w} NP_{sk}$ for all $s$ and $w$ (32)

$CNP_{sW} \geq MinFraction_s * CN_{sW}$ for all $s$ **(optional)** (33)

Constraints 34-36 relate to penalty 6. Constraints 34-35 use $CNP_{sw}$ to compute $Diff_{sw}$ which is summed in (36) to give $TimingPenalty_s$, a measure of the irregularity in the timing of section $s$'s mass meetings. The right sides of (34-35) are opposites, and $Diff_{sw}$ must exceed both, so the minimization of $TimingPenalty_s$ in (1) means that $Diff_{sw}$ equals the absolute value of the difference between the terms on the right of (34) at an optimal solution. Term 1 on the right of (34) is the actual number of mass meetings for section $s$ in weeks 1 to $w$ combined ($CNP_{sw}$), and term 2 is the prorated number of such meetings. For example, consider a section that normally meets once a week during an 8-week semester but has only four mass meetings in weeks 5, 6, 7, and 8 during a pandemic. In this case, $CNP_{sw}$ is (0, 0, 0, 0, 1, 2, 3, 4) and the prorated number of meetings is (0.5, 1, 1.5, 2, 2.5, 3, 3.5, 4) in weeks 1 to 8 respectively. The absolute value of the difference is (0.5, 1, 1.5, 2, 1.5, 1, 0.5, 0), and these terms sum to 8 = $TimingPenalty_s$. This high penalty indicates the meetings are poorly distributed. If the mass meetings occur in weeks 2, 4, 6, and 8, $CNP_{sw}$ is (0, 1, 1, 2, 2, 3, 3, 4), and the absolute value of the difference between it and the prorated number of meetings is (0.5, 0, 0.5, 0, 0.5, 0, 0.5, 0) which sums to 2. This low penalty indicates the meetings are well distributed. The expected value of $TimingPenalty_s$ grows linearly in the number of mass meetings for section $s$, so the factor $CN_{sW}$ is omitted from objective component 6 (see Eq. 1).

$Diff_{sw} \geq CNP_{sw} - \left(\frac{w - FirstW_s + 1}{LastW_s - FirstW_s + 1}\right) CNP_{sW}$ for all $s$ and all $w$ from $FirstW_s$ to $LastW_s$ (34)

$Diff_{sw} \geq \left(\frac{w - FirstW_s + 1}{LastW_s - FirstW_s + 1}\right) CNP_{sW} - CNP_{sw}$ for all $s$ and all $w$ from $FirstW_s$ to $LastW_s$ (35)

$TimingPenalty_s = \sum_{w=FirstW_s}^{LastW_s} Diff_{sw}$ for all $s$ (36)

Constraint 37 relates to objective penalty 7. It forces $U_s$ to be 0 if $CNP_{sW} < (MinFraction_s * CN_{sW})$, i.e., if less than $MinFraction_s$ of section $s$'s meetings are in person (note the capital $W$).

$CNP_{sW} \geq U_s * (MinFraction_s * CN_{sW})$ for all $s$ (37)

## 6. OFFICE algorithm

The math model is prohibitively large for real problem instances. For example, our case study has 172 rooms, 1834 sections, 590 weekly meeting times, and 16 weeks. This amounts to approximately 3 billion (=172*1834*590*16) $Z_{rsmw}$ variables. Preliminary experiments with IBM ILOG CPLEX 12.9 and only 18 course sections crashed our computer on 100% of occasions due to a lack of memory. This motivated the development of a heuristic algorithm for the HP-UCSP.

Our heuristic algorithm is called OFFICE (Optimized Face-to-Face Interaction for Commerce and Education). It is designed to generate excellent course schedules that balance the seven competing objectives in (1). Nevertheless, the algorithm has priorities. Roughly speaking, the algorithm prioritizes objective components 1-7 in the following order: 7, 2, 3, 1, 4, 5, 6. Fairness is paramount in our approach, so the algorithm's main goal is to allow a basic fraction of meetings for every course section to be in person (objective 7). However, room assignment possibilities are limited a-priori to ensure good results for objectives 1-4. Once objective component 7 is zeroed out, OFFICE tries to maximize total student-hours in a classroom across the entire campus (objective 5). Built into its core scheduling engine is a preference to schedule meetings in weeks that are spread throughout the semester (objective 6).

The OFFICE algorithm follows the seven-step procedure shown in Table 7. In step 1, we create an exhaustive list of *possible room assignments* (PRAs) for each section *s*. Each PRA is a unique set of one or more rooms with sufficient capacity to accommodate the section, no wasted rooms, and acceptable values for *NumRooms*$_s$, *DistPenalty*$_s$, *PrefPenalty*$_s$, and *WastedSeats*$_s$ (see Table 5).  In step 2, the PRAs for each section are ranked according to a weighted average of parts 1-4 of the objective function.

**Table 7.** OFFICE algorithm procedure

| Step | Description |
|---|---|
| **1** | Form an exhaustive list of possible room assignments (PRAs) for each section *s*. Each PRA for section *s* is a unique set of rooms that: <br> • Has enough total capacity to host a mass meeting for section *s* <br> • Uses no more than *MaxRooms* rooms and has no wasted rooms <br> • Has a value of $DistPenalty_s$ no greater than $MaxDistPenalty_s$ <br> • Has a value of $PrefPenalty_s$ no greater than $MaxPrefPenalty_s$ <br> • Has a value of $WastedSeats_s$ no greater than $MaxWastedSeats_s$ |
| **2** | Rank the PRAs for section *s* according to a weighted average of parts 1-4 of the objective function. |
| **3** | Form a **master list** of all sections' in-person **meetings**. Each meeting in the list is defined by a section number and **default week** when the meeting occurs during a normal semester. |
| **4** | Assign a **wave number** from 1 to *V* to each meeting in the master list. The lower the number, the higher the priority for scheduling. The first *V*-1 waves together have $[MinFraction_s * CN_{sW}]$ meetings for each section *s*. Wave *V* has all other meetings, i.e., $CN_{sW} - [MinFraction_s * CN_{sW}]$ meetings for each section *s*. |
| **5** | Begin with a blank **master schedule** in which no mass meetings are scheduled. Let *v* = 1. |
| **6** | Place the meetings in wave *v* (*v* < *V*) into the schedule as follows: <br> **A** Without disturbing the mass meetings scheduled in previous waves, call the **scheduling engine** (Table 8) to create a **current solution** in which as many wave *v* meetings as possible are placed in the schedule, assuming the best PRA (with the lowest penalty) for each wave *v* section is used. <br> **B** Compute the objective value of the current solution according to (1). <br> **C** If *TL* seconds have elapsed, STOP and display the best feasible solution found. <br> **D** If all wave *v* meetings are in the schedule and *v* < *V*-1 (*v* = *V*-1), increase *v* by 1 and go to step 6 (7). <br> **E** Create a **next solution** by (i) removing all mass meetings for DU(*low*, *high*) wave *v* sections from the current solution, (ii) selecting a new, random PRA for each section just removed, and (iii) calling the **scheduling engine** to place as many unscheduled wave *v* meetings as possible in the schedule. <br> **F** Compute the objective value difference between the next and current solutions. <br> **G** Use simulated annealing to decide if the next solution replaces the current solution. Go to step 6B. |
| **7** | Place the wave *V* meetings in the schedule using the procedure described in steps 6A-6G except that step 6D, and the PRA phrase in step 6A, do not apply. After *TL* sec have elapsed, STOP. Display the best soln. found. |

Mass meetings are placed into the *master schedule* (i.e., schedule) in *V* waves (steps 3-7) with more important and harder-to-schedule meetings assigned to earlier waves.  Each meeting belonging to wave *v* must be feasibly placed in the schedule as a mass meeting before any meeting in wave *v*+1 is considered.  At no time may the number of mass meetings in the schedule for any section be less than the number in the schedule at the end of any previous wave.  In each wave, every feasible solution is a master schedule for the entire semester that is evaluated according to (1), and simulated annealing principles decide if neighboring solutions are accepted.  For each section *s*, $[MinFraction_s * CN_{sW}]$ of its meetings belong to one of the first *V*-1 waves, and the remaining $CN_{sW} - [MinFraction_s * CN_{sW}]$ of its meetings are in wave *V*.  If all meetings in waves 1 to *V*-1 are feasibly placed in the schedule, fairness is achieved and the value of objective component 7 is 0.  Thereafter, in wave *V*, the algorithm places as many remaining meetings as possible into the schedule (objective 5) while searching for better room assignments (objectives 1-4) and meeting distributions across the weeks of the semester (objective 6).

A key feature of the OFFICE algorithm is its procedure of adding items to the schedule one mass meeting at a time and removing them one section at a time.  The first feasible solution in each wave is created by calling a *scheduling engine* which adds mass meetings to the last feasible solution from the previous wave, one meeting at a time (step 6A). Thereafter, each neighboring solution is created by (i) removing all mass meetings for DU(*low*, *high*) sections from the schedule; (ii) selecting a new, random PRA for each section just removed; and then (iii) calling the *scheduling engine* to schedule as many unscheduled meetings as possible (step 6E).  The algorithm terminates after *TL* seconds have elapsed.

Table 8 shows the scheduling engine procedure. The first step is to create a randomly scrambled list of all meetings in wave *v* that are not in the master schedule. Each such meeting is defined by a section number and *default week* when it occurs during a normal semester. We then proceed through the list and attempt to schedule each meeting, one at a time, during its default week and in rooms specified by the PRA most recently selected for its section. All eligible meeting times in that week (see $I_{sm}$ in Table 4) are considered. If any meeting time works, the meeting is added to the schedule. After all meetings in the list are considered, we revisit those not yet in the schedule, and we try to add each of them, one at a time, to the schedule one week before or after the meeting's default week. After considering all meetings in this manner, we revisit those still not in the schedule and try to place each of them, one at a time, in the schedule two weeks before or two weeks after the default week. This continues until all meetings in the list are scheduled (i.e., wave *v* is completed) or no more meetings can be feasibly scheduled in any week.

**Table 8.** Scheduling engine procedure

| Step | Description |
|---|---|
| 1 | Create a randomly scrambled list of all meetings in wave $v$ that are not yet in the master schedule. Let $Sect_n$ and $Week_n$ be the section number and default week for meeting $n$ in the list. Let $Deviation = 0$. |
| 2 | Proceed through the list and try to schedule each meeting $n$, one at a time, in either week $Week_n + Deviation$ or $Week_n - Deviation$, and in the rooms corresponding to the current PRA being considered for section $Sect_n$. All eligible meeting times $m$ (given by $I_{sm}$ where $s = Sect_n$) in the week(s) at hand are considered. |
| 3 | If any meeting time is available for meeting $n$, place a mass meeting in the master schedule during that time and delete the meeting from the list. If weeks $Week_n + Deviation$ and $Week_n - Deviation$ are both out of bounds, delete meeting $n$ from the list. |
| 4 | If any meetings remain in the list, increase $Deviation$ by 1 and go to step 2. If not, STOP. |

## 7. Case study: [*UniversityXYZ*] campus in Fall 2022

Our first two consultations with the [*UniversityXYZ*] Registrar's Office took place in late February and early March 2020, just weeks before the COVID-19 outbreak. Back then we were interested in regular classroom assignment, but our interests later expanded to include pandemic preparedness and response. In this section, we describe our work to develop schedules for the fall 2022 semester at [*UniversityXYZ*]. We discuss our dataset, describe our experimental setup, and present the results of two experiments that consider course scheduling (i) in normal times and (ii) during a pandemic.

*7.1. Data collection*

[*UniversityXYZ*] is [*StateXYZ*]'s second largest university with 23,000 students, 206 academic programs, and a 104-acre campus. It is one of two "R1" universities in the state. Highly competent staff members in the registrar's office gave us data for the fall 2022 semester one week before it began.

Table 9 summarizes the data. The campus has 172 general classrooms, and the equipment in each room—a podium computer with internet access and an LCD projector—is suitable for any of the 1834 general course sections. Each section is *nonredundant*—cross-listed sections with the same meeting days, time, and instructor are already combined—and at least one in-person meeting is planned for each section. Ten large chemistry courses have frequent daytime lectures (e.g., MWF from 11:30-12:20) and infrequent in-person evening exams in a different room (e.g., from 5:30-7 PM on Sept. 26, Oct. 24, Nov. 14, and Dec. 12). Given the disparity in meeting frequency, location, time, and duration between the lectures and exams, we created two sections for each such course: one for the lectures and one for the exams.

Instruction takes place six days a week over 16 weeks with minimal instruction in the 16[th] week and on Saturdays. All meeting times start at :00 or :30 after the hour, so 30-minute timeslots are used. The classrooms reside in 17 buildings, all within 800 meters of each other. The total seating capacity of all rooms during normal times (a pandemic) is 9953 (2531). Ninety-seven organizations offer courses, and their preferences for teaching in different buildings are categorized into four *tiers* 1, 2, 3, 4 with buildings in tier 1 being most desirable. The average number of buildings in each organization's tier (1, 2, 3, 4) is about (1, 2, 3, 11). Buildings in tier 4 should be avoided when possible, so we decided to use building preference penalty values ($P_{gb}$) of (0, 1, 2, 6) for tiers (1, 2, 3, 4) respectively (Table 9, bottom).

**Table 9.** Summary of the in-person instruction planned for fall 2022 at [*UniversityXYZ*]

| General data | | | Notes | | | |
|---|---|---|---|---|---|---|
| S | 1834 course sections | | 1782 (52) regular (irregular) sections. 42,827 total in-person meetings. | | | |
| R | 172 classrooms | | 9953 (2531) total usable seats during normal times (a pandemic) | | | |
| M | 590 weekly meeting times | | Equipment in all rooms is compatible with all sections (all $J_{rs} = 1$) | | | |
| W | 16 weeks | | Week 16 is only used by ESL (English as 2nd language) sections | | | |
| D | 96 academic days | | 7 academic days are holidays | | | |
| T | 156 timeslots per week (each 30 min) | | 1 room is not available during first 16 days due to flooding | | | |
| B | 17 buildings | | Min (max) distance between two buildings is 58.2 (783.7) meters | | | |
| G | 97 organizations offering courses | | 1,714,290 total student-hours of in-person attendance in fall 2022 | | | |

| Section data | | | | | | |
|---|---|---|---|---|---|---|
| Course level (catalog number) | No. sections | Enrollment | No. sections | No. in-person meetings in a normal semester | No. sections (* = all irregular) | |
| 0 (062-099) | 62 | 0 | 12 | 1-12 | 41* | |
| 1 (100-199) | 656 | 1-19 | 790 | 23, 32, or 33 | 3* | |
| 2 (200-299) | 334 | 20-39 | 787 | 13-15 (1 mtg / week) | 753 | |
| 3 (300-399) | 316 | 40-59 | 97 | 27-31 (2 mtgs / week) | 901 | |
| 4 (400-499) | 176 | 60-99 | 75 | 41-42 (3 mtgs / week) | 62 | |
| 5-6 (500-699) | 144 | 100-149 | 60 | 55-61 (4 mtgs / week) | 65 | |
| 7-9 (700-999) | 146 | 150+ | 13 | 69 (5 mtgs / week) | 9 | |

| Meeting day | No. sections | Meeting start time | No. sections | Meeting duration | No. sections |
|---|---|---|---|---|---|
| Monday | 731 | 08:00-09:30 | 275 | 31 - 60 min (2 timeslots) | 601 |
| Tuesday | 686 | 10:00-11:30 | 592 | 61 - 90 min (3 timeslots) | 803 |
| Wednesday | 762 | 12:00-13:30 | 307 | 91-120 min (4 timeslots) | 149 |
| Thursday | 680 | 14:00-15:30 | 301 | 121-150 min (5 timeslots) | 3 |
| Friday | 235 | 16:00-17:30 | 315 | 151-180 min (6 timeslots) | 267 |
| Saturday | 5 | 18:00-19:30 | 44 | 181+ min (7+ timeslots) | 11 |

| Room and building data | | | | | |
|---|---|---|---|---|---|
| Normal capacity | No. rooms | Covid capacity | No. rooms | | |
| 16-29 | 64 | 5-9 | 76 | | |
| 30-39 | 40 | 10-19 | 67 | | |
| 40-49 | 29 | 20-29 | 11 | | |
| 50-99 | 15 | 30-39 | 9 | | |
| 100-199 | 16 | 40-70 | 9 | | |
| 200+ | 8 | | | | |

| Building preference penalty ($P_{gb}$) | No. g-b pairs | No. rooms in the building | No. buildings | Floor where the room is located | No. rooms |
|---|---|---|---|---|---|
| | | 2 | 4 | 0 (basement) | 27 |
| 0 | 112 | 3-5 | 4 | 1 | 94 |
| 1 | 188 | 6-10 | 3 | 2 | 34 |
| 2 | 266 | 11-20 | 4 | 3 | 16 |
| 6 | 1083 | 21+ | 2 | 5 | 1 |

*7.2. General experimental setup*

Data received from the registrar was translated into text files, and the OFFICE algorithm was coded into MS Visual C++ 2010. All experiments were run on a desktop PC with the Windows 7 Enterprise operating system, an Intel core i7-4770 (3.40 GHz) processor, and 16 GB of RAM.

We set up the OFFICE algorithm with $V = 5$ waves. Wave 1 has all meetings for each of the ten sections representing midterm exams for large chemistry courses (Section 7.1). ($MinFraction_s = 1.0$ for these special sections.) Wave 2 has $MinFraction_s$ of the meetings for each of the (42 other) sections with irregular in-person meetings. Most of these sections have infrequent in-person meetings, so opportunities for mass meetings are scarce. Wave 3 (4) has $MinFraction_s$ of the meetings for each regular section that has an enrollment of at least (below) $ENR$, and wave 5 has all other meetings.

Table 10 shows the parameter values used in the experiments. These values were identified during preliminary experiments not described here. Parameters in the math model are at the top; those used in the OFFICE algorithm are at the bottom. Table 11 describes in simple terms how the objective function changes when a section's mass meetings change. For example, teaching a section in one more building (100) has twice the penalty as moving it from a tier 1 to tier 2 building (50). The high values of $\alpha_5$ and $\alpha_7$ in Table 10 reflect our goal of scheduling at least $MinFraction_s$ of every section's meetings in person, and the low start temperature for pandemic scheduling means that inferior neighbors are almost never accepted. The values of $ImportanceOfLevel_l$ show that lower-level (e.g., first year) sections are prioritized over higher-level sections. For regular classroom assignment, the values (1, 6, 465) for $MaxRooms$, $MaxPrefPenalty_s$, and $MaxWastedSeats_s$ respectively mean that every room with enough seats for a section is a feasible PRA for the section. For pandemic scheduling, the values of $MaxDistPenalty_s$ and $MaxPrefPenalty_s$—480 and 2 respectively for most sections—are the lowest values that allow at least one feasible PRA to be identified for nearly all sections. These values are higher than (480, 2) for (3, 13) sections, respectively, to avoid infeasibility. The values of *low* and *high* indicate that 1-10 sections' mass meetings are removed from the schedule each iteration (step 6E in Table 7).

**Table 10.** Parameter values used in the experiments

| Math model parameter | Value for regular classroom assignment | Value for scheduling during a pandemic |
|---|---|---|
| $\alpha_1$ | * | 15 |
| $\alpha_2$ ($\alpha_{21}, \alpha_{22}, \alpha_{23}, \alpha_{24}, \alpha_{25}$) | * | 1 (1, 100, 10, 30, 3) |
| $\alpha_3$ | 50 | 50 |
| $\alpha_4$ | 0 | 0 |
| $\alpha_5$ (*Exp*)** | 1000 (1) | 1000 (1) |
| $\alpha_6$ | 100 | 100 |
| $\alpha_7$ | 1,000,000 | 1,000,000 |
| $ImportanceOfLevel_l$ | (4, 5, 4, 3, 2, 1) for level (0, 1, 2, 3, 4, 5-9) | (4, 5, 4, 3, 2, 1) for level (0, 1, 2, 3, 4, 5-9) |
| $MaxRooms$ | 1 | 5 |
| $MaxDistPenalty_s$ | * | (480, 700, 1000) for (1831, 1, 2) sections |
| $MaxPrefPenalty_s$ | 6 | (2, 6) for (1821, 13) sections*** |
| $MaxWastedSeats_s$ | 465 | 20**** |
| $MinFraction_s$ | 1.0 for all sections | all sections equal (except wave 1 sects. = 1.0) |
| OFFICE algorithm parameter | | |
| $V$ | 5 | 5 |
| (*low, high*) | (1, 10) | (1, 10) |
| $ENR$ | 80 | 80 |
| Start temperature | 20,000,000 (reset to 200 at start of wave 5) | 200 |
| Temperature factor | .99999 | .999999 |
| $TL$ | 3600 seconds | 3600 seconds |

\* Values are irrelevant because only one room is assigned to each section
\*\* Values shown are used in wave 5 of the algorithm only. In waves 1-4, the value of $\alpha_5$ (*Exp*) is 40,000 (3).
\*\*\* Value is changed to 6 if a room is already assigned by the registrar and it has a building preference penalty of 6
\*\*\*\* Value is changed to 465 if a room is already assigned by the registrar and it wastes more than 20 seats

**Table 11.** Impact of changing a section's mass meetings

| Change | Relative amount added to objective value |
| --- | --- |
| One more room is used | 15 |
| Distance between buildings used increases by 1 meter | 1 |
| One more building is used | 100 |
| Total "floor distance" increases by 1 | 10 |
| One more floor is used | 30 |
| Room assignment has one less pair of adjacent rooms | 3 |
| Building preference penalty increases by 1 | 50 |
| Room assignment has one more wasted seat | 0 |
| 10% more meetings are online instead of in person | 100 |
| Timing of one mass meeting is off by 1 more week | 4 (approximate) |
| Section no longer has $MinFraction_s$ of its meetings in person | 1,000,000 |

*7.3. Experiment 1: Regular classroom assignment*

Table 12 shows the results of six independent, one-hour runs of the OFFICE algorithm for regular classroom assignment (task 4 in Figure 1). The second column shows the performance of the registrar's room assignments that were used in fall 2022. Column 3 shows the performance of OFFICE's room assignments, averaged over the six runs. Column 4 shows the results for OFFICE in the best of the six runs. Objective values (in thousands) are at the top, and other results are at the bottom. The values of objective components 5 and 7 are 0 in all cases, indicating that 100% of all sections' meetings are feasibly scheduled. Components 1-2 are 0 because every section meets in one room, and component 6 is nonzero owing to the discrete, imperfect way of computing $TimingPenalty_s$. Regarding objective component 3, the average value achieved by OFFICE (238,730) is less than half the value achieved in the registrar's assignment (484,832). Overall, OFFICE finds significantly better room assignments than the registrar's software program.

**Table 12.** Experiment 1 results

| Value (in thousands) of… | Registrar's room assignments | OFFICE algorithm (Average of 6 runs) | OFFICE algorithm (Best of 6 runs) |
|---|---|---|---|
| Objective components 1, 2, 4, 5, 7 | 0 | 0 | 0 |
| Objective component 3 | 484,832 | 238,730 | 201,761 |
| Objective component 6 | 343,542 | 343,542 | 343,542 |
| Total objective value | 828,374 | 582,272 | 545,303 |
| No. sections with building preference penalty 0 | 1205 | 1363 | 1386 |
| No. sections with building preference penalty 1 | 378 | 346 | 343 |
| No. sections with building preference penalty 2 | 156 | 94 | 79 |
| No. sections with building preference penalty 6 | 95 | 30 | 26 |
| Avg. no. PRAs generated per section | 1.0 | 138.1 | 138.1 |

*7.4. Experiment 2: Pandemic preparedness and response*

We now consider pandemic preparedness and response in two *scenarios*: (R) if each section's room assignment must include the room initially assigned to it by the registrar (task 6 in Fig. 1) and (NR) if not (task 4A in Fig. 1). For each scenario, we explore six values of *MinFraction$_s$* (for the course sections in waves 2-4)—0.05, 0.1, 0.15, 0.2, 0.25, and 0.3—and we run OFFICE six times for each value, yielding 72 (1-hour) runs altogether. In all runs, *MinFraction$_s$* is 1.0 for all course sections in wave 1.

The main result of this experiment is that, in all runs with *MinFraction$_s$* ≤ 0.25, OFFICE was able to schedule at least *MinFraction$_s$* of all sections' meetings in person. Figure 3 shows the average objective value (in thousands) of the solutions identified by OFFICE in these runs. Here, performance improves as the decision maker faces fewer restrictions. For example, OFFICE finds better solutions in scenario NR (gray line with round points) than R (black line with square points). Also, the objective value improves as *MinFraction$_s$* decreases, i.e., as the conditions for zeroing out objective 7 are relaxed.

The results when *MinFraction$_s$* = 0.3 are not shown in Figure 3 because no such runs resulted in all sections reaching the 30% threshold. The best (average) of six runs for scenario R still left 93 (144.0) sections with less than 30% of their meetings in person, and the best (average) of six runs for scenario NR still left 1 (2.5) section(s) below the 30% threshold. Although no run resulted in all sections reaching the 30% threshold, the results for scenario NR are clearly better than scenario R.

**Figure 3.** OFFICE algorithm performance in experiment 2 when *MinFraction$_s$* ranges from 0.05 to 0.25 (averaged over 6 runs). In all runs, at least *MinFraction$_s$* of all sections' meetings are in person.

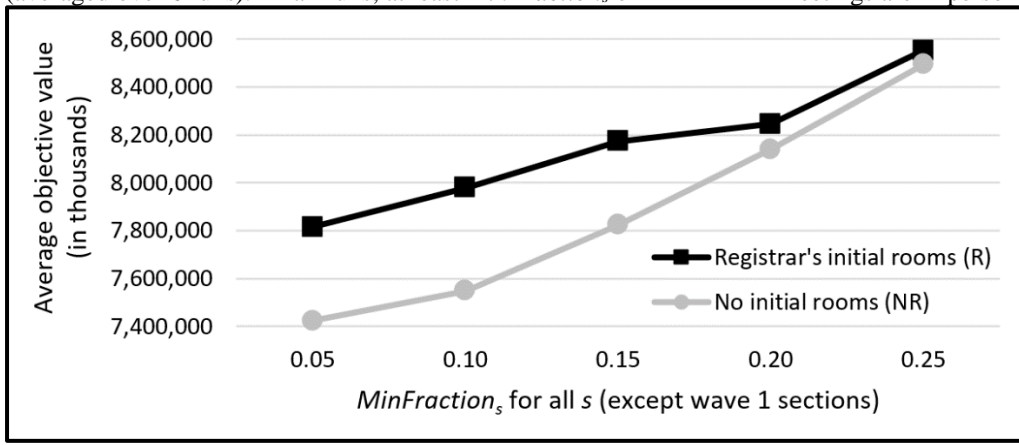

Table 13 shows the best solutions obtained by OFFICE when $MinFraction_s = 0.25$ for scenario R (left) and NR (right). These high-quality solutions show that OFFICE is an effective scheduling tool that successfully balances numerous objectives. The results are better for scenario NR because the decision maker has more room assignment options (814.9 per section on average) than in scenario R (417.0).

Regarding objective components 1-2, the R solution assigns 2.14 rooms, that occupy 1.22 floors and 1.04 buildings, to the average section whereas the NR solution assigns 2.04 rooms, that occupy 1.15 floors and 1.03 buildings, to the average section. If we weight each section by the number of students enrolled, we find that the average student attending a section in the R solution (NR solution) sees 2.64 rooms that occupy 1.35 floors and 1.16 buildings (2.62 rooms that occupy 1.31 floors and 1.15 buildings).

Regarding objective 3, the R solution assigns more sections to nonpreferred buildings than the NR solution, but it performs well given that it must incorporate the registrar's initial room assignments. Indeed, the registrar assigns (378, 156, 95) sections to penalty (1, 2, 6) buildings respectively (Table 12), and the R solution assigns only a few more sections respectively: (401, 166, 95). Importantly, in both solutions most sections have favorable building assignments despite the large number of undesirable organization-building pairs in the dataset (bottom of Table 9). Regarding objective 6, the R solution has slightly better meeting timing than the NR solution (*TimingPenalty* = 380,745 vs. 390,354), but both solutions perform well given that *TimingPenalty* = 343,542 in a perfectly scheduled semester (Table 12).

**Table 13.** Best solutions obtained by OFFICE when $MinFraction_s = 0.25$

|  | Registrar's initial rooms (R) (Best of 6 runs) | No initial rooms (NR) (Best of 6 runs) |
|---|---|---|
| Avg. no. PRAs generated per section | 417.0 | 814.9 |
| **Value (in thousands) of…** | | |
| Objective component 1 | 321,362 | 316,682 |
| Objective component 2 | 649,934 | 577,250 |
| Objective component 3 | 557,745 | 458,634 |
| Objective component 4 | 0 | 0 |
| Objective component 5 | 6,546,222 | 6,527,042 |
| Objective component 6 | 380,745 | 390,354 |
| Objective component 7 | 0 | 0 |
| Total objective value | 8,456,008 | 8,269,962 |
| No. rooms per section (section avg, student avg) | (2.14, 2.64) | (2.04, 2.62) |
| No. floors per section (section avg, student avg) | (1.22, 1.35) | (1.15, 1.31) |
| No. buildings per section (section avg, student avg) | (1.04, 1.16) | (1.03, 1.15) |
| No. sections with building preference penalty 0 | 1172 | 1224 |
| No. sections with building preference penalty 1 | 401 | 357 |
| No. sections with building preference penalty 2 | 166 | 227 |
| No. sections with building preference penalty 6 | 95 | 26 |
| Total student-hours of mass meeting attendance | 844,736.5 (49.3%) | 847,585.5 (49.4%) |
| % of (small, large) sections with <30% of meetings in person | (16.9%, 43.8%) | (17.3%, 58.9%) |
| % of (small, large) sections with ≥80% of meetings in person | (24.9%, 9.6%) | (27.1%, 11.0%) |

In both solutions, more than 840,000, or 49%, of the 1,714,290 student-hours of in-person attendance in the original schedule can still occur in a classroom if a pandemic occurs. For comparison, if each section remains in the room assigned to it by the registrar and rotating attendance is used, the total amount of in-person attendance—including the residential spread, hybrid split, and hybrid touch point delivery modes identified by Navabi-Shirazi et al. (2022)—is 770,832.5 student-hours, including 127,000 student-hours of residential spread (i.e., mass meeting) attendance. Thus, OFFICE generates a nearly sevenfold increase in mass meeting attendance, and a 10% increase in total in-person attendance, compared to a good alternative. Plus, it ensures that at least 25% of each section's meetings are in person.

The final rows of Table 13 show the disparity in the fraction of meetings held in person across the sections. In the R solution, 16.9% (24.9%) of the small sections—with fewer than 100 students—have less than 30% (at least 80%) of their meetings in person. On the other hand, 43.8% (9.6%) of the large sections—with at least 100 students—have less than 30% (at least 80%) of their meetings in person. The NR solution has even greater disparities. These results show the importance of maintaining a 25% baseline for all sections to prevent the wholesale sacrifice of large sections for the benefit of smaller ones.

## 8. Conclusion

Scheduling university courses is especially difficult when social distancing requirements reduce the useful capacity of classrooms during a pandemic. In this work, we introduce a new framework for university course scheduling that is based on three principles which have received little attention in the literature: fairness, simultaneous attendance, and high-precision scheduling. Our approach is fair by ensuring that a certain fraction (e.g., 25%) of the instruction in every course section occurs in the classroom during a pandemic. Our focus on simultaneous attendance, in which all students in a section meet at the same time in one or more nearby rooms, creates opportunities for in-person midterm exams and group activities that are not available with rotating attendance. As discussed earlier, the best way to achieve fairness with simultaneous attendance is to schedule at high precision across all days of a semester rather than a single, repeating week. These ideas led to the development of the OFFICE algorithm and its testing on a real case with 1834 in-person course sections, 172 classrooms, and 96 days. Results show the algorithm is an effective course scheduling tool in both normal and pandemic times.

Future work might apply the OFFICE algorithm to other institutions or to event/space management in crowded office buildings where the demand for space exceeds the supply. A detailed comparison of our approach to a method that uses rotating attendance and room reassignments would also be worthwhile. Future research may also explore the new course delivery modes introduced in Section 3.

**References**

Barnhart, C., Bertsimas, D., Delarue, A., & Yan, J. (2022). Course scheduling under sudden scarcity: Applications to pandemic planning. *Manufacturing & Service Operations Management*, 24(2), 727-745.

Bettinelli, A., Cacchiani, V., Roberti, R., & Toth, P. (2015). An overview of curriculum-based course timetabling. *TOP*, 23(2), 313-349.

Freeman, S., Eddy, S. L., McDonough, M., Smith, M. K., Okoroafor, N., Jordt, H., & Wenderoth, M. P. (2014). Active learning increases student performance in science, engineering, and mathematics. *Proceedings of the National Academy of Sciences*, 111(23), 8410-8415.

Gore, A. B., Kurz, M. E., Saltzman, M. J., Splitter, B., Bridges, W. C., & Calkin, N. J. (2022). Clemson University's rotational attendance plan during COVID-19. *INFORMS Journal on Applied Analytics*, 52(6), 553-567.

Johnes, J. (2015). Operational research in education. *European J. of Operational Research*, 243, 683-696.

Johnson, C., & Wilson, R. L. (2022). Practice summary: A multiobjective assignment model for optimal socially distanced classrooms for the Spears School of Business at Oklahoma State University. *INFORMS Journal on Applied Analytics*, 52(3), 295-300.

Khamechian, M., & Petering, M. E. H. (2022). A mathematical modeling approach to university course planning. *Computers & Industrial Engineering*, 168, article number 107855.

Lindahl, M., Mason, A. J., Stidsen, T., & Sørensen, M. (2018). A strategic view of university timetabling. *European Journal of Operational Research*, 266, 35-45.

McCollum, B., McMullan, P., Parkes, A. J., Burke, E. K., & Qu, R. (2012). A new model for automated examination timetabling. *Annals of Operations Research*, 194, 291-315.

McCollum, B., Schaerf, A., Paechter, B., McMullan, P., Lewis, R., Parkes, A. J., Gaspero, L. D., Qu, R., & Burke, E. K. (2010). Setting the research agenda in automated timetabling: The second international timetabling competition. *INFORMS Journal on Computing*, 22(1), 120-130.

Moug, K., Padharia, M., Smith, K., Shen, S., Denton, B., & Cohn, A. (2022). Mid-semester pandemic-driven course rescheduling with integer programming. *Proceedings of the IISE Annual Conference*, 2022.

Navabi-Shirazi, M., El Tonbari, M., Boland, N., Nazzal, D., & Steimle, L. N. (2022). Multicriteria course mode selection and classroom assignment under sudden space scarcity. *Manufacturing & Service Operations Management*, 24(6), 2797-3306.

Phillips, A. E., Waterer, H., Ehrgott, M., & Ryan, D. M. (2015). Integer programming methods for large-scale practical classroom assignment problems. *Computers & Operations Research*, 53, 42-53.

Pillay, N. (2014). A survey of school timetabling research. *Annals of Operations Research*, 218, 261-293.

Reeves, P. (2020). Can social distancing drive better classroom utilization policies in higher ed? https://www.gordian.com/resources/social-distancing-and-classroom-utilization, website accessed May 10, 2023.

Rudová, H., Müller, T., & Murray, K. (2011). Complex university course timetabling. *Journal of Scheduling*, 14, 187-207.

WHO (World Health Organization) (2023). WHO coronavirus (COVID-19) dashboard. https://covid19.who.int, website accessed May 10, 2023.

Schady, N., Holla, A., Sabarwal, S., Silva, J., & Chang, A. Y. (2023). Collapse and recovery: How the COVID-19 pandemic eroded human capital and what to do about it, *The World Bank*.